\documentclass[10pt,a4paper]{amsart}
\usepackage{amsmath}
\usepackage{amssymb}
\usepackage{amsthm}
\usepackage{amsfonts}
\usepackage{microtype}
\usepackage{mathrsfs}
\usepackage{graphicx}
\usepackage{color}
\usepackage{latexsym}
\usepackage{rotating}
\usepackage{xspace}
\usepackage[all]{xy}
\usepackage{longtable}
\usepackage{leftidx}
\usepackage{mathtools}
\usepackage{tikz}
\usetikzlibrary{arrows}
\usetikzlibrary{decorations.markings}
\usepackage{geometry}

\newtheorem{theorem}{Theorem}[section]
\numberwithin{equation}{theorem}
\newtheorem{lemma}[theorem]{Lemma}

\newtheorem{corollary}[theorem]{Corollary}

\theoremstyle{definition}
\newtheorem{definition}[theorem]{Definition}
\newtheorem{example}[theorem]{Example}

\theoremstyle{conjecture}

\newtheorem{question}[theorem]{Question}

\newcommand{\im}{\operatorname{im}}

\newcommand{\ara}{\operatorname{ara}}

\newcommand{\cd}{\operatorname{cd}}

\newcommand{\hd}{\operatorname{hd}}

\newcommand{\V}{\operatorname{V}}

\newcommand{\Ext}{\operatorname{Ext}}
\newcommand{\Supp}{\operatorname{Supp}}

\newcommand{\Tor}{\operatorname{Tor}}
\newcommand{\Hom}{\operatorname{Hom}}

\newcommand{\Ann}{\operatorname{Ann}}

\newcommand{\fa}{\frak{a}}

\newcommand{\suchthat}{\;\ifnum\currentgrouptype=16 \middle\fi|\;}

\newenvironment{prf}[1][Proof]{\begin{proof}[\bf #1]}{\end{proof}}

\begin{document}

\author[K. Divaani-Aazar, H. Faridian and M. Tousi]{Kamran Divaani-Aazar, Hossein Faridian
and Massoud Tousi}

\title[Local homology, finiteness of Tor ...]
{Local homology, finiteness of Tor modules and cofiniteness}

\address{K. Divaani-Aazar, Department of Mathematics, Alzahra University, Vanak, Tehran, Iran, Zip Code 19834; and
School of Mathematics, Institute for Research in Fundamental Sciences (IPM), Tehran, Iran, P.O. Box 19395-5746.}
\email{kdivaani@ipm.ir}

\address{H. Faridian, Department of Mathematics, Shahid Beheshti University, G. C., Evin, Tehran, Iran, Zip Code 1983963113.}
\email{h.faridian@yahoo.com}

\address{M. Tousi, Department of Mathematics, Shahid Beheshti University, G. C., Evin, Tehran, Iran, Zip Code 1983963113.}
\email{mtousi@ipm.ir}


\keywords {Cofinite modules; local cohomology; local homology.\\
The research of the first and third authors are supported by grants from IPM (No. 94130212 and No. 92130211; respectively).}

\begin{abstract}
Let $\frak a$ be an ideal of a commutative noetherian ring $R$ with unity and $M$ an $R$-module supported at $\V(\fa)$. Let $n$
be the supermum of the integers $i$ for which $H^{\fa}_i(M)\neq 0$. We show that $M$ is $\fa$-cofinite if and only if the $R$-module
$\Tor^R_i(R/\fa,M)$ is finitely generated for every $0\leq i\leq n$. This provides a hands-on and computable finitely-many-steps
criterion to examine $\mathfrak{a}$-confiniteness. Our approach relies heavily on the theory of local homology which
demonstrates the effectiveness and indispensability of this tool.
\end{abstract}

\maketitle

\section{Introduction}

Throughout this note, $R$ denotes a commutative noetherian ring with unity.

The theory of local cohomology has been developed tremendously during the six decades of investigation after its introduction by Grothendieck. However, its dual theory, i.e. the
theory of local homology has not been developed much. The theory of local homology was initiated by Matlis \cite{Ma} in 1974, and its study was continued by Simon in \cite{Si1} and \cite{Si2}. After the groundbreaking works of Greenlees and May \cite{GM}, and Alonso Tarr\'{i}o, Jerem\'{i}as
L\'{o}pez and Lipman \cite {AJL}, a new era in the study of this theory has begun; see e.g. \cite{Sc}, \cite{CN1}, \cite{CN2},
\cite{Fr} and \cite{Ri}.

To prove Theorem \ref{1.5} below, we deploy the theory of local homology as the main tool, which in effect illustrates the conspicuous role
of this theory in the study of cofiniteness, and provides along the way results of independent interest for a deeper insight into the
theory of local homology.

It is proved in \cite[Lemma 2.3]{BA} that if $\mathfrak{a}$ is an ideal of $R$ that can be generated by $n$ elements and $M$ is an $R$-module,
then the $R$-module $\Tor^{R}_{i}(R/\mathfrak{a},M)$ is finitely generated for every $i \geq 0$ if and only if it is finitely generated
for every $0 \leq i \leq n$. Recall that the cohomological dimension of $M$ with respect to $\mathfrak{a}$ is $$\cd(\mathfrak{a},M):=\sup
\{i\geq 0 \suchthat H^{i}_{\mathfrak{a}}(M)\neq 0 \},$$ and note that $\cd(\mathfrak{a},M) \leq \ara(\fa)$. From this angle, the following
generalization of this result is posed as a question in \cite[Question 2]{PAB}.

\begin{question} \label{1.2}
Let $\mathfrak{a}$ be an ideal of $R$. If the $R$-module $\Tor^{R}_{i}(R/\mathfrak{a},M)$ is finitely generated for every $0 \leq i \leq
\cd(\mathfrak{a},R)$, then is the $R$-module $\Tor^{R}_{i}(R/\mathfrak{a},M)$ finitely generated for every $i \geq 0$?
\end{question}

We answer this question in the affirmative in an even stronger form. Given an ideal $\mathfrak{a}$ of $R$ and an $R$-module $M$, we define
the homological dimension of $M$ with respect to $\fa$ by $$\hd(\mathfrak{a},M):= \sup \{i\geq 0 \suchthat H^{\mathfrak{a}}_{i}(M)\neq 0 \}.$$
It follows from \cite[Theorem 2.5 and Corollary 3.2]{GM}, that $\hd(\fa,M)\leq \cd(\fa,R)$. As a matter of fact, we derive the following result.

\begin{theorem} \label{1.5}(See Theorem \ref{2.9}.)
Let $\mathfrak{a}$ be an ideal of $R$ and $M$ an $R$-module. Then the following assertions are equivalent:
\begin{enumerate}
\item[(i)] $H^{\mathfrak{a}}_{i}(M)$ is a finitely generated $\widehat{R}^{\mathfrak{a}}$-module for every $i\geq 0$.
\item[(ii)] $\Tor^{R}_{i}(R/\mathfrak{a},M)$ is a finitely generated $R$-module for every $i\geq 0$.
\item[(iii)] $\Tor^{R}_{i}(R/\mathfrak{a},M)$ is a finitely generated $R$-module for every $0\leq i\leq \hd(\mathfrak{a},M)$.
\end{enumerate}
\end{theorem}

In his algebraic geometry seminar of 1962, Grothendieck asked whether the $R$-modules $\Hom_{R}\left(R/\mathfrak{a},
H^{i}_{\mathfrak{a}}(M)\right)$ were finitely generated for every ideal $\mathfrak{a}$ of $R$ and every finitely generated
$R$-module $M$; see  \cite[Conjectures 1.1 and 1.2]{Gr}. Two years later, Hartshorne provided a counterexample in
\cite[Section 3]{Ha}, to show that this question does not have an affirmative answer in general. He then proposed the
following definition.

\begin{definition} \label{1.1}
Let $\mathfrak{a}$ be an ideal of $R$. An $R$-module $M$ is said to be $\mathfrak{a}$-\textit{cofinite} if $\Supp_{R}(M)
\subseteq \V(\mathfrak{a})$, and $\Ext^{i}_{R}(R/\mathfrak{a},M)$ is finitely generated for every $i \geq 0$.
\end{definition}

As an immediate application of Theorem \ref{1.5} , we deduce the following result. It is worth noting that Tor modules are easier to compute and
more handy to deal with than local homology modules. Moreover, the next result sets forth a finitely-many-steps criterion for
cofiniteness.

\begin{corollary} \label{1.6}(See Corollary \ref{2.11}.)
Let $\mathfrak{a}$ be an ideal of $R$ and $M$ an $R$-module with $\Supp_{R}(M) \subseteq \V(\mathfrak{a})$. Then the following assertions
are equivalent:
\begin{enumerate}
\item[(i)] $M$ is $\mathfrak{a}$-cofinite.
\item[(ii)] $\Tor^{R}_{i}(R/\mathfrak{a},M)$ is a finitely generated $R$-module for every $0 \leq i \leq \hd(\mathfrak{a},M)$.
\end{enumerate}
\end{corollary}

\section{Main Results}

We need to equip ourselves with a handful of lemmas and propositions in advance. We begin with recalling the definition of local homology modules.
Given an ideal $\fa$ of $R$, we let $\Lambda^{\mathfrak{a}}(M):=\widehat{M}^{\mathfrak{a}}=\varprojlim (M/\mathfrak{a}^{n}M)$ for an $R$-module $M$,
and $\Lambda^{\mathfrak{a}}(f):=\widehat{f}$ for an $R$-homomorphism $f:M\rightarrow N$. This provides us with the so-called $\mathfrak{a}$-adic completion functor
$\Lambda^{\mathfrak{a}}(-)$ on the category of $R$-modules. Moreover, we define the $i$th local homology module of an $R$-module $M$ with respect to $\mathfrak{a}$ to be
$H^{\mathfrak{a}}_{i}(M):= H_{i}\left(\Lambda^{\mathfrak{a}}(P)\right)$ for every $i\geq 0$, where $P$ is any projective resolution of $M$.

\begin{lemma} \label{2.1}
Let $\mathfrak{a}$ be an ideal of $R$ and $M$ an $R$-module. If $M/\mathfrak{a}M$ is a finitely generated $R$-module, then $\widehat{M}^{\mathfrak{a}}$
is a finitely generated $\widehat{R}^{\mathfrak{a}}$-module.
\end{lemma}

\begin{prf}
As $M/\mathfrak{a}M$ is finitely generated, there is a finitely generated submodule $N$ of $M$ such that $M=N+\mathfrak{a}M$. Let $\iota:N\rightarrow M$
be the inclusion map. Then by \cite[Lemma 1.2]{Si1}, the $\widehat{R}^{\mathfrak{a}}$-homomorphism $\widehat{\iota}:\widehat{N}^{\mathfrak{a}}\rightarrow
\widehat{M}^{\mathfrak{a}}$ is surjective.
On the other hand, since $\widehat{N}^{\mathfrak{a}} \cong \widehat{R}^{\mathfrak{a}} \otimes_{R}N$, it is clear that $\widehat{N}^{\mathfrak{a}}$ is
a finitely generated $\widehat{R}^{\mathfrak{a}}$-module. It then follows that $\widehat{M}^{\mathfrak{a}}$ is a finitely generated
$\widehat{R}^{\mathfrak{a}}$-module.
\end{prf}

We recall that an $R$-module $M$ is said to be $\mathfrak{a}$-adically quasi-complete if the completion map
$\theta_{M}:M\rightarrow \widehat{M}^{\mathfrak{a}}$ is an epimorphism.

\begin{lemma} \label{2.2}
Let $\mathfrak{a}$ be an ideal of $R$. Then the following assertions hold:
\begin{enumerate}
\item[(i)] Any submodule of an $\mathfrak{a}$-adically separated $R$-module is $\mathfrak{a}$-adically separated.
\item[(ii)] Any homomorphic image of an $\mathfrak{a}$-adically quasi-complete $R$-module is $\mathfrak{a}$-adically quasi-complete.
\item[(iii)] If $f:M\rightarrow N$ is a homomorphism of $\mathfrak{a}$-adically complete $R$-modules, then both $\ker f$ and $\im f$ are
$\mathfrak{a}$-adically complete.
\end{enumerate}
\end{lemma}

\begin{prf}
(i): As $\ker \theta_{M}=\bigcap_{i=1}^{\infty} \mathfrak{a}^{i}M$ for any $R$-module $M$, the assertion is clear.

(ii): If $M$ is $\mathfrak{a}$-adically quasi-complete and $f:M\rightarrow N$ is an epimorphism, then it follows from \cite[Lemma 1.2]{Si1} that $\widehat{f}:\widehat{M}^{\mathfrak{a}}\rightarrow \widehat{N}^{\mathfrak{a}}$ is surjective. Therefore, the commutative diagram
\[\begin{tikzpicture}[every node/.style={midway},]
  \matrix[column sep={3em}, row sep={2.5em}]
  {\node(1) {$M$}; & \node(2) {$N$}; \\
  \node(3) {$\widehat{M}^{\mathfrak{a}}$}; & \node(4) {$\widehat{N}^{\mathfrak{a}}$}; \\};
  \draw[decoration={markings,mark=at position 1 with {\arrow[scale=1.5]{>}}},postaction={decorate},shorten >=0.5pt] (1) -- (2) node[anchor=south] {$f$};
  \draw[decoration={markings,mark=at position 1 with {\arrow[scale=1.5]{>}}},postaction={decorate},shorten >=0.5pt] (3) -- (4) node[anchor=south] {$\widehat{f}$};
  \draw[decoration={markings,mark=at position 1 with {\arrow[scale=1.5]{>}}},postaction={decorate},shorten >=0.5pt] (1) -- (3) node[anchor=west] {$\theta_{M}$};
  \draw[decoration={markings,mark=at position 1 with {\arrow[scale=1.5]{>}}},postaction={decorate},shorten >=0.5pt] (2) -- (4) node[anchor=west] {$\theta_{N}$};
\end{tikzpicture}\]
shows that $\theta_{N}$ is surjective, i.e. $N$ is $\mathfrak{a}$-adically quasi-complete.

(iii): As $N$ is $\mathfrak{a}$-adically complete, it can be seen by inspection that
$$\ker f= \bigcap_{i=1}^{\infty} \left(\ker f + \mathfrak{a}^{i}M\right),$$
i.e. $\ker f$ is a closed submodule of $M$ in the $\mathfrak{a}$-adic topology. It now follows from \cite[Proposition 1.3 (ii)]{Si1} that $\ker f$
is $\mathfrak{a}$-adically complete. On the other hand, $\im f$ is both a submodule of $N$ and a homomorphic image of $M$, so using (i) and (ii),
we infer that $\im f$ is $\mathfrak{a}$-adically complete.
\end{prf}

\begin{lemma} \label{2.3}
Let $\mathfrak{a}$ be an ideal of $R$ and $M$ an $R$-module. Suppose that $R$ is $\mathfrak{a}$-adically complete, $M$ is $\mathfrak{a}$-adically
quasi-complete, and $M/\mathfrak{a}M$ is a finitely generated $R$-module. Let $\theta_{M}:M \rightarrow \widehat{M}^{\mathfrak{a}}$ be the
completion map with $K= \ker \theta_{M}$. Then $\widehat{K}^{\mathfrak{a}}=0=H^{\mathfrak{a}}_{0}(K)$.
\end{lemma}

\begin{prf}
Since $M/\mathfrak{a}M$ is a finitely generated $R$-module, Lemma \ref{2.1} implies that $\widehat{M}^{\mathfrak{a}}$ is a finitely generated
$R$-module. From the short exact sequence
$$0\rightarrow K \rightarrow M \xrightarrow {\theta_{M}} \widehat{M}^{\mathfrak{a}} \rightarrow 0,$$
we get the exact sequence
$$\Tor^{R}_{1}\left(R/\mathfrak{a},\widehat{M}^{\mathfrak{a}}\right)\rightarrow K/\mathfrak{a}K \rightarrow M/\mathfrak{a}M,$$
which implies that $K/\mathfrak{a}K$ is a finitely generated $R$-module. Subsequently, from the short exact sequence
$$0\rightarrow K/\mathfrak{a}K \rightarrow M/\mathfrak{a}K \rightarrow \widehat{M}^{\mathfrak{a}} \rightarrow 0,$$
we deduce that $M/\mathfrak{a}K$ is a finitely generated $R$-module. It then follows that the zero submodule of $M/\mathfrak{a}K$ has a minimal
primary decomposition
$$\mathfrak{a}K/\mathfrak{a}K= \bigcap_{i=1}^{n}Q_{i}/\mathfrak{a}K,$$
which in turn gives a minimal primary decomposition
$$\mathfrak{a}K= \bigcap_{i=1}^{n}Q_{i}$$
of $\mathfrak{a}K$.
We prove that $K=\mathfrak{a}K$, by showing that $K \subseteq Q_{i}$ for every $1 \leq i \leq n$. Assume to the contrary that there is an integer
$1 \leq j \leq n$ and an element $x\in K \backslash Q_{j}$. Then we have $\mathfrak{a}x \subseteq \mathfrak{a}K \subseteq Q_{j}$. This means that
the homomorphism $M/Q_{j} \xrightarrow {a} M/Q_{j}$ is not injective for every $a\in \mathfrak{a}$. As $Q_{j}$ is a primary submodule of $M$, we
conclude that $\mathfrak{a} \subseteq \sqrt{(Q_{j}:_{R}M)}$. Therefore, there is an integer $t \geq 1$ such that $\mathfrak{a}^{t}M \subseteq Q_{j}$.
But $$K=\bigcap_{i=1}^{\infty}\mathfrak{a}^{i}M \subseteq \mathfrak{a}^{t}M \subseteq Q_{j},$$
which is a contradiction. Therefore, $K=\mathfrak{a}K$, and so by \cite[Lemma 5.1 (ii)]{Si1}, $\widehat{K}^{\mathfrak{a}}=0=H^{\mathfrak{a}}_{0}(K)$.
\end{prf}

The next result may be of independent interest.

\begin{lemma} \label{2.4}
Let $\mathfrak{a}$ be an ideal of $R$ and $M$ an $R$-module. Then for any $j\geq 0$ we have:
\[
    H_{i}^{\mathfrak{a}}\left(H_{j}^{\mathfrak{a}}(M)\right)\cong
\begin{dcases}
    H_{j}^{\mathfrak{a}}(M) & \text{if } i= 0\\
    0              & \text{if } i \geq 1.
\end{dcases}
\]
\end{lemma}

\begin{prf}
Let
$$F= \cdots \rightarrow F_{2} \xrightarrow {\partial^{F}_{2}} F_{1} \xrightarrow {\partial^{F}_{1}} F_{0} \rightarrow 0$$
be a free resolution of $M$. Then by definition $$H^{\mathfrak{a}}_{j}(M)= \ker \widehat{\partial^{F}_{j}}/ \im  \widehat{
\partial^{F}_{j+1}},$$for any $j\geq 0$. By Lemma \ref{2.2} (iii), both $\im \widehat{\partial^{F}_{j+1}}$ and $\ker
\widehat{\partial^{F}_{j}}$ are $\mathfrak{a}$-adically complete. Therefore, invoking \cite[Lemmas 5.1 (i), and 5.2 (i)]{Si1},
there are natural isomorphisms $$\alpha: H_{0}^{\mathfrak{a}}\left(\im \widehat{\partial^{F}_{j+1}}\right) \rightarrow \im
\widehat{\partial^{F}_{j+1}},$$and$$\beta: H_{0}^{\mathfrak{a}}\left(\ker \widehat{\partial^{F}_{j}}\right) \rightarrow \ker
\widehat{\partial^{F}_{j}}.$$ Now the short exact sequence
\begin{equation} \label{eq:2.4.1}
0 \rightarrow  \im  \widehat{\partial^{F}_{j+1}} \xrightarrow {\iota} \ker \widehat{\partial^{F}_{j}} \xrightarrow {\pi}
H^{\mathfrak{a}}_{j}(M) \rightarrow 0,
\end{equation}
yields the following commutative diagram with exact rows
\[\begin{tikzpicture}[every node/.style={midway},]
  \matrix[column sep={3em}, row sep={2.5em}]
  {\node(1) {}; & \node(2) {$H_{0}^{\mathfrak{a}}\left(\im  \widehat{\partial^{F}_{j+1}}\right)$}; & \node(3) {$H_{0}^{\mathfrak{a}}
  \left(\ker  \widehat{\partial^{F}_{j}}\right)$}; & \node(4) {$H_{0}^{\mathfrak{a}}\left(H^{\mathfrak{a}}_{j}(M)\right)$}; &
  \node(5) {$0$}; \\
  \node(6) {$0$}; & \node(7) {$\im  \widehat{\partial^{F}_{j+1}}$}; & \node(8) {$\ker \widehat{\partial^{F}_{j}}$}; & \node(9)
  {$H^{\mathfrak{a}}_{j}(M)$}; & \node(10) {$0$}; \\};
  \draw[decoration={markings,mark=at position 1 with {\arrow[scale=1.5]{>}}},postaction={decorate},shorten >=0.5pt] (2) -- (3)
  node[anchor=south] {$H_{0}^{\mathfrak{a}}(\iota)$};
  \draw[decoration={markings,mark=at position 1 with {\arrow[scale=1.5]{>}}},postaction={decorate},shorten >=0.5pt] (3) -- (4)
  node[anchor=south] {$H_{0}^{\mathfrak{a}}(\pi)$};
  \draw[decoration={markings,mark=at position 1 with {\arrow[scale=1.5]{>}}},postaction={decorate},shorten >=0.5pt] (4) -- (5)
  node[anchor=south] {};
  \draw[decoration={markings,mark=at position 1 with {\arrow[scale=1.5]{>}}},postaction={decorate},shorten >=0.5pt] (6) -- (7)
  node[anchor=south] {};
  \draw[decoration={markings,mark=at position 1 with {\arrow[scale=1.5]{>}}},postaction={decorate},shorten >=0.5pt] (7) -- (8)
  node[anchor=south] {$\iota$};
  \draw[decoration={markings,mark=at position 1 with {\arrow[scale=1.5]{>}}},postaction={decorate},shorten >=0.5pt] (8) -- (9)
  node[anchor=south] {$\pi$};
  \draw[decoration={markings,mark=at position 1 with {\arrow[scale=1.5]{>}}},postaction={decorate},shorten >=0.5pt] (9) -- (10)
  node[anchor=south] {};
  \draw[decoration={markings,mark=at position 1 with {\arrow[scale=1.5]{>}}},postaction={decorate},shorten >=0.5pt] (2) -- (7)
  node[anchor=west] {$\cong$} node[anchor=east] {$\alpha$};
  \draw[decoration={markings,mark=at position 1 with {\arrow[scale=1.5]{>}}},postaction={decorate},shorten >=0.5pt] (3) -- (8)
  node[anchor=west] {$\cong$} node[anchor=east] {$\beta$};
  \draw[dashed,decoration={markings,mark=at position 1 with {\arrow[scale=1.5]{>}}},postaction={decorate},shorten >=0.5pt] (4) --
  (9) node[anchor=west] {$\gamma$} node[anchor=east] {$\exists$};
\end{tikzpicture}\]
from which we deduce that $H_{0}^{\mathfrak{a}}(\iota)$ is injective and $\gamma$ is an isomorphism.
On the other hand, the short exact sequence ~\eqref{eq:2.4.1} yields an exact sequence
$$0= H_{i+1}^{\mathfrak{a}}\left(\ker \widehat{\partial^{F}_{j}}\right) \rightarrow H_{i+1}^{\mathfrak{a}}
\left(H^{\mathfrak{a}}_{j}(M)\right) \rightarrow H_{i}^{\mathfrak{a}}\left(\im  \widehat{\partial^{F}_{j+1}}\right) =0,$$
for every $i \geq 1$, where the vanishing follows from \cite[Lemma 5.2 (i)]{Si1}.
It follows that $H_{i+1}^{\mathfrak{a}}\left(H^{\mathfrak{a}}_{j}(M)\right)=0$ for every $i \geq 1$. It further yields the exact
sequence
$$0= H_{1}^{\mathfrak{a}}\left(\ker \widehat{\partial^{F}_{j}}\right) \rightarrow H_{1}^{\mathfrak{a}}\left(H^{\mathfrak{a}}_{j}(M)
\right) \rightarrow H_{0}^{\mathfrak{a}}\left(\im  \widehat{\partial^{F}_{j+1}}\right) \xrightarrow {H_{0}^{\mathfrak{a}}(\iota)}
H_{0}^{\mathfrak{a}}\left(\ker \widehat{\partial^{F}_{j}}\right).$$
As $H_{0}^{\mathfrak{a}}(\iota)$ is injective, we conclude that $H_{1}^{\mathfrak{a}}\left(H^{\mathfrak{a}}_{j}(M)\right)=0$.
\end{prf}

\begin{lemma} \label{2.5}
Let $\mathfrak{a}$ be an ideal of $R$, $M$ an $R$-module, and $s\geq 0$ an integer. Suppose that $R$ is $\mathfrak{a}$-adically complete, and
$(R/\mathfrak{a})\otimes_{R}H^{\mathfrak{a}}_{s}(M)$ is a finitely generated $R$-module. Then $H^{\mathfrak{a}}_{s}(M)$ is a finitely
generated $R$-module.
\end{lemma}

\begin{prf}
Let $L=H^{\mathfrak{a}}_{s}(M)$. As $L/\mathfrak{a}L$ is finitely generated, Lemma \ref{2.1} implies that $\widehat{L}^{\mathfrak{a}}$
is a finitely generated $R$-module. Lemma \ref{2.2} yields that $L$ is $\mathfrak{a}$-adically quasi-complete, so Lemma \ref{2.3} implies
that $H_{0}^{\mathfrak{a}}(\ker \theta_{L})=0$, where $\theta_{L}:L \rightarrow \widehat{L}^{\mathfrak{a}}$
is the completion map. Hence from the short exact sequence
$$0 \rightarrow \ker \theta_{L} \rightarrow L \xrightarrow {\theta_{L}} \widehat{L}^{\mathfrak{a}} \rightarrow 0,$$
we get the exact sequence
$$0=H_{0}^{\mathfrak{a}}(\ker \theta_{L}) \rightarrow H_{0}^{\mathfrak{a}}(L) \rightarrow H_{0}^{\mathfrak{a}}\left(\widehat{L}^
{\mathfrak{a}}\right) \rightarrow 0,$$
implying that $H_{0}^{\mathfrak{a}}(L)\cong H_{0}^{\mathfrak{a}}\left(\widehat{L}^{\mathfrak{a}}\right)$. By \cite[Lemma 5.2 (i)]{Si1},
we have $H_{0}^{\mathfrak{a}}\left(\widehat{L}^{\mathfrak{a}}\right)\cong \widehat{L}^{\mathfrak{a}}$. On the other hand, Lemma \ref{2.4}
implies that $H_{0}^{\mathfrak{a}}(L)\cong L$. Hence $L\cong \widehat{L}^{\mathfrak{a}}$ is a finitely generated $R$-module.
\end{prf}

\begin{lemma} \label{2.6}
Let $\mathfrak{a}$ be an ideal of $R$ and $M$ an $R$-module. Then there is a convergent spectral sequence
$$E^{2}_{p,q}= \Tor^{R}_{p}\left(R/\mathfrak{a},H^{\mathfrak{a}}_{q}(M)\right) \underset {p} \Rightarrow \Tor^{R}_{p+q}(R/\mathfrak{a},M).$$
\end{lemma}

\begin{prf}
Let $\mathcal{F}=(R/\mathfrak{a})\otimes_{R}-$, and $\mathcal{G}=\Lambda^{\mathfrak{a}}(-)$. Then $\mathcal{F}$ is right exact, and
$\mathcal{G}(P)$ is left $\mathcal{F}$-acyclic for every projective $R$-module $P,$ since the $\mathfrak{a}$-adic completion of a flat
$R$-module is flat by \cite[1.4.7]{B}. Therefore, by \cite[Theorem 10.48]{Ro}, there is a Grothendieck spectral sequence
$$E^{2}_{p,q}= L_{p}\mathcal{F}\left(L_{q}\mathcal{G}(M)\right) \underset {p} \Rightarrow L_{p+q}(\mathcal{FG})(M).$$
Let $F$ be a free resolution of $M$. By \cite[Theorem 1.3.1]{B} or \cite[Theorem 15]{Ma}, we have
$$(R/\mathfrak{a})\otimes_{R}\Lambda^{\mathfrak{a}}(F)\cong (R/\mathfrak{a})\otimes_{R}F,$$
whence
$$L_{p+q}(\mathcal{FG})(M) = H_{p+q}\left((\mathcal{FG})(F)\right) \cong \Tor^{R}_{p+q}(R/\mathfrak{a},M).$$
\end{prf}

\begin{lemma} \label{2.7}
Let $\mathfrak{a}$ be an ideal of $R$, $M$ an $R$-module, and $s \geq 0$ an integer. If $\Tor^{R}_{i}\left(R/\mathfrak{a},H^{\mathfrak{a}}_{j}(M)\right)$
is finitely generated for every $i \geq 0$ and $0 \leq j \leq s-1$, and $\Tor^{R}_{s}(R/\mathfrak{a},M)$ is finitely generated, then
$(R/\mathfrak{a})\otimes_{R} H^{\mathfrak{a}}_{s}(M)$ is finitely generated.
\end{lemma}

\begin{prf}
By Lemma \ref{2.6}, there is a convergent spectral sequence
$$E^{2}_{p,q}= \Tor^{R}_{p}\left(R/\mathfrak{a},H^{\mathfrak{a}}_{q}(M)\right) \underset {p} \Rightarrow \Tor^{R}_{p+q}(R/\mathfrak{a},M).$$
Therefore, there is a finite filtration
$$0=U^{-1} \subseteq U^{0} \subseteq \cdots \subseteq U^{s} = \Tor^{R}_{s}(R/\mathfrak{a},M),$$
such that $U^{p}/U^{p-1}\cong E^{\infty}_{p,s-p}$ for every $p\in \mathbb{Z}$. As $\Tor^{R}_{s}(R/\mathfrak{a},M)$ is finitely generated,
we conclude that $$E^{\infty}_{0,s}\cong U^{0}/U^{-1} = U^{0}$$ is finitely generated.
Let $r \geq 2$, and consider the differentials
$$E^{r}_{r,s-r+1}\xrightarrow {d^{r}_{r,s-r+1}} E^{r}_{0,s}\xrightarrow {d^{r}_{0,s}} E^{r}_{-r,s+r-1}=0,$$
where the vanishing comes from the facts that $E^{2}_{-r,s+r-1}=0$ and $E^{r}_{-r,s+r-1}$ is a subquotient of $E^{2}_{-r,s+r-1}$ for every
$r \geq 2$. On the other hand, as $s-r+1 \leq s-1$, the hypothesis implies that $E^{r}_{r,s-r+1}$ and consequently $\im d^{r}_{r,s-r+1}$
are finitely generated for every $r\geq 2$.
We thus obtain
$$E^{r+1}_{0,s} \cong \ker d^{r}_{0,s}/ \im d^{r}_{r,s-r+1} = E^{r}_{0,s}/ \im d^{r}_{r,s-r+1},$$
and consequently a short exact sequence
\begin{equation} \label{eq:2.7.1}
0 \rightarrow \im d^{r}_{r,s-r+1} \rightarrow E^{r}_{0,s} \rightarrow E^{r+1}_{0,s} \rightarrow 0.
\end{equation}
There is an integer $r_{0}\geq 2$, such that $E^{\infty}_{0,s}=E^{r+1}_{0,s}$ for every $r\geq r_{0}$. It follows that $E^{r_{0}+1}_{0,s}$
is finitely generated.
Now the short exact sequence ~\eqref{eq:2.7.1} implies that $E^{r_{0}}_{0,s}$ is finitely generated. Using the short exact sequence
~\eqref{eq:2.7.1} inductively, we conclude that $E^{2}_{0,s}=(R/\mathfrak{a})\otimes_{R} H^{\mathfrak{a}}_{s}(M)$ is finitely generated as desired.
\end{prf}

The next result is of grave importance in the theory of cofinite modules, which is included here for the convenience of the reader.

\begin{lemma} \label{2.8}
Let $\mathfrak{a}$ be an ideal of $R$ and $M$ an $R$-module. Then the following conditions are equivalent:
\begin{enumerate}
\item[(i)] $\Ext_R^i(R/\mathfrak{a},M)$ is a finitely generated $R$-module for every $i\geq 0$.
\item[(ii)]  $\Ext_R^i(N,M)$ is a finitely generated $R$-module for every finitely generated $R$-module $N$ with $\Supp_{R}(N)\subseteq
 \V(\mathfrak{a})$ and every $i\geq 0$.
\item[(iii)] $\Tor^R_i(R/\mathfrak{a},M)$ is a finitely generated $R$-module for every $i\geq 0$.
\item[(iv)] $\Tor^R_i(N,M)$ is a finitely generated $R$-module for every finitely generated $R$-module $N$ with $\Supp_{R}(N)\subseteq
\V(\mathfrak{a})$ and every $i\geq 0$.
\item[(v)] $H^{\mathfrak{a}}_i(M)$ is a finitely generated $\widehat{R}^{\mathfrak{a}}$-module for every $i\geq 0$.
\end{enumerate}
\end{lemma}

\begin{prf}
Apply \cite[Propositions 7.1, 7.2, and Theorem 7.4]{WW} to the module case, bearing in mind that
$H^{\mathfrak{a}}_{i}(M)=0$ for every $i> \cd(\mathfrak{a},R)$.
\end{prf}

Now, we are ready to prove our main result.

\begin{theorem} \label{2.9}
Let $\mathfrak{a}$ be an ideal of $R$ and $M$ an $R$-module. Then the following assertions are equivalent:
\begin{enumerate}
\item[(i)] $H^{\mathfrak{a}}_{i}(M)$ is a finitely generated $\widehat{R}^{\mathfrak{a}}$-module for every $i\geq 0$.
\item[(ii)] $\Tor^{R}_{i}(R/\mathfrak{a},M)$ is a finitely generated $R$-module for every $i\geq 0$.
\item[(iii)] $\Tor^{R}_{i}(R/\mathfrak{a},M)$ is a finitely generated $R$-module for every $0\leq i\leq \hd(\mathfrak{a},M)$.
\end{enumerate}
\end{theorem}

\begin{prf}
(i) $\Rightarrow$ (ii): Holds by Lemma \ref{2.8}.

(ii) $\Rightarrow$ (iii): Obvious.

(iii) $\Rightarrow$ (i): We claim that we may assume that $R$ is $\mathfrak{a}$-adically complete. Indeed, if $\Tor^{R}_{i}(R/
\mathfrak{a},M)$ is a finitely generated $R$-module,
then by \cite[Lemma 2.3]{Si2}, $$\Tor^{\widehat{R}^{\mathfrak{a}}}_{i}\left(\widehat{R}^{\mathfrak{a}}/\mathfrak{a}\widehat{R}^
{\mathfrak{a}},\widehat{R}^{\mathfrak{a}}\otimes_{R}M\right) \cong \Tor^{R}_{i}(R/\mathfrak{a},M)$$
is a finitely generated $R$-module, and thus a finitely generated $\widehat{R}^{\mathfrak{a}}$-module. If the result is proved in
the complete case, then by \cite[Lemma 2.3]{Si2},
$$H^{\mathfrak{a}}_{i}(M) \cong H^{\mathfrak{a}\widehat{R}^{\mathfrak{a}}}_{i}\left(\widehat{R}^{\mathfrak{a}}\otimes_{R}M\right)$$
is a finitely generated $\widehat{R}^{\mathfrak{a}}$-module. In addition, the latter isomorphism indicates that $$\hd(\fa,M)=\hd\left(\fa \widehat{R}^{\mathfrak{a}},\widehat{R}^{\mathfrak{a}}\otimes_{R}M\right).$$ Having assumed that $R$ is $\mathfrak{a}$-adically complete,
we proceed to do induction on $s$ to show that if $\Tor^{R}_{i}(R/\mathfrak{a},M)$ is a finitely generated $R$-module for every
$0 \leq i \leq s$,
then $H^{\mathfrak{a}}_{s}(M)$ is a finitely generated $R$-module. Clearly, this establishes the implication.
If $s=0$, then $M/\mathfrak{a}M$ is finitely generated. Now \cite[Theorem 1.1 and Lemma 5.1 (iv)]{Si1} implies that
$$\begin{array}{ll} H_{0}^{\mathfrak{a}}(M)/\mathfrak{a}H_{0}^{\mathfrak{a}}(M)&\cong \widehat{H_{0}^{\mathfrak{a}}(M)}^{\mathfrak{a}}/\mathfrak{a}
\widehat{H_{0}^{\mathfrak{a}}(M)}^{\mathfrak{a}}\\
&\cong  \widehat{M}^{\mathfrak{a}}/\mathfrak{a}\widehat{M}^{\mathfrak{a}}\\
&\cong  M/\mathfrak{a}M.
\end{array}$$
Therefore, the $R$-module $H_{0}^{\mathfrak{a}}(M)/\mathfrak{a}H_{0}^{\mathfrak{a}}(M)$ is finitely generated. Hence Lemma
\ref{2.5} implies that $H_{0}^{\mathfrak{a}}(M)$ is finitely generated. Now assume that $s> 0$ and the claim is true for $s-1$.
Hence $H^{\mathfrak{a}}_{i}(M)$ is finitely generated for every $0\leq i\leq s-1$. It follows that $\Tor^{R}_{i}\left(R/\mathfrak{a},
H^{\mathfrak{a}}_{j}(M)\right)$ is finitely generated for every $i\geq 0$ and $0 \leq j \leq s-1$, and $\Tor^{R}_{s}(R/\mathfrak{a},M)$
is finitely generated. Lemma \ref{2.7} implies that $(R/\mathfrak{a})\otimes_{R}H_{s}^{\mathfrak{a}}(M)$ is finitely generated, and
subsequently Lemma \ref{2.5} yields that $H_{s}^{\mathfrak{a}}(M)$ is finitely generated as desired.
\end{prf}

\begin{corollary} \label{2.10}
Let $\mathfrak{a}$ be an ideal of $R$ and $M$ an $R$-module. If the $R$-module $$\bigoplus_{i=0}^{\hd(\mathfrak{a},R)}\Tor^{R}_{i}(R/
\mathfrak{a},M)$$ is finitely generated, then $\Tor^{R}_{i}(N,M)$ is a finitely generated $R$-module for every finitely generated
$R$-module $N$ with $\Supp_{R}(N) \subseteq \V(\mathfrak{a})$ and every $i \geq 0$.
\end{corollary}

\begin{prf} The proof is clear in light of Theorem \ref{2.9} and Lemma \ref{2.8}.
\end{prf}

The next result provides an applicable criterion for cofiniteness.

\begin{corollary} \label{2.11}
Let $\mathfrak{a}$ be an ideal of $R$ and $M$ an $R$-module with $\Supp_{R}(M) \subseteq \V(\mathfrak{a})$. Then the following assertions
are equivalent:
\begin{enumerate}
\item[(i)] $M$ is $\mathfrak{a}$-cofinite.
\item[(ii)] $\Tor^{R}_{i}(R/\mathfrak{a},M)$ is a finitely generated $R$-module for every $0 \leq i \leq \hd(\mathfrak{a},M)$.
\end{enumerate}
\end{corollary}

\begin{prf} The proof is clear in view of Theorem \ref{2.9} and Lemma \ref{2.8}.
\end{prf}

\begin{corollary} \label{2.12}
Let $N$ be a finitely generated $R$-module and  $M$ an $R$-module. If $\Tor^{R}_{i}\left(R/\Ann_RN,M\right)$ is a finitely generated
$R$-module for every $0 \leq i \leq \hd\left(\Ann_RN,M\right)$, then $\Tor^{R}_{i}(N,M)$ is a finitely generated $R$-module for every $i\geq 0$.
\end{corollary}

\begin{prf} The proof is clear in light of Theorem \ref{2.9} and Lemma \ref{2.8}.
\end{prf}



\begin{thebibliography}{99}

\bibitem[AJL]{AJL}{L. Alonso Tarr\'{i}o, A. Jerem\'{i}as L\'{o}pez and J. Lipman}, {\it Local homology and cohomology on schemes}, Ann.
Sci. \'{E}cole Norm. Sup., (4){\bf 30}(1), (1997), 1-39.

\bibitem[B]{B}{J. Bartijn}, {\it Flatness, completion, regular sequences, un m\'{e}nage \`{a} trois}, Thesis. Utrecht (1985).

\bibitem[BA]{BA}{K. Bahmanpour and M. Aghapournahr}, {\it A note on cofinite modules},  Comm. Algebra, \textbf{44}(9), (2016), 3683-3691.

\bibitem[CN1]{CN1}{N.T. Cuong and T.T. Nam}, {\it A local homology theory for linearly compact modules}, J. Algebra,  {\bf 319}(11), (2008),
4712-4737.

\bibitem[CN2]{CN2}{N.T. Cuong and T.T. Nam}, {\it The I-adic completion and local homology for Artinian modules}, Math. Proc. Cambridge Philos. Soc.,
{\bf 131}(1), (2001), 61-72.

\bibitem[Fr]{Fr}{A. Frankild}, {\it Vanishing of local homology}, Math. Z., {\bf 244}(3), (2003), 615-630.

\bibitem[Gr]{Gr}{A. Grothendieck}, {\it Cohomologie locale des faisceaux coh\`{e}rents et th\`{e}or\`{e}mes de Lefschetz locaux
et globaux}, (SGA 2), North-Holland, Amsterdam, 1968.

\bibitem[GM]{GM}{J.P.C. Greenlees and J.P. May}, {\it Derived functors of $I$-adic completion and local homology}, J. Algebra,
\textbf{149}(2), (1992), 438-453.

\bibitem[Ha]{Ha}{R. Hartshorne}, {\it Affine duality and cofiniteness}, Invent. Math., \textbf{9} (1969/1970), 145-164.

\bibitem[Ma]{Ma}{E. Matlis}, {\it The higher properties of $R$-sequences}, J. Algebra, \textbf{50}(1), (1978), 77-112.

\bibitem[PAB]{PAB}{G. Pirmohammadi, K. Ahmadi Amoli and K. Bahmanpour}, {\it Cofiniteness with respect to ideals of cohomological
dimension at most one}, preprint.

\bibitem[Ri]{Ri}{A.S. Richardson}, {\it Co-localization, co-support and local homology}, Rocky Mountain J. Math., {\bf 36}(5), (2006), 1679-1703.

\bibitem[Ro]{Ro}{J.J. Rotman}, {\it An introduction to homological algebra}, Universitext. Springer, New York, second edition, 2009.

\bibitem[Si1]{Si1}{A.M. Simon}, {\it Some homological properties of complete modules}, Math. Proc. Camb. Phil. Soc., \textbf{108}(2),
(1990), 231-246.

\bibitem[Si2]{Si2}{A.M. Simon}, {\it Adic-completion and some dual homological results}, Publ. Mat. Camb., \textbf{36}(2B), (1992),
965-979.

\bibitem[WW]{WW}{S. Sather-Wagstaff and R. Wicklein}, {\it Support and adic finiteness for complexes}, Arxive: 1401.6925v3, [math.AC]
5 Jun. 2015.

\bibitem[Sc]{Sc}{P. Schenzel}, {\it Proregular sequences, local cohomology, and completion}, Math. Scand., {\bf 92}(2), (2003), 161-180.

\end{thebibliography}
\end{document}